\documentclass[leqno]{macroicaiado}
\DeclareSymbolFont{rsfscript}{OMS}{rsfs}{m}{n}
\DeclareSymbolFontAlphabet{\mathrsfs}{rsfscript}
 \DeclareMathOperator{\ad}{ad}
\DeclareMathOperator{\chro}{\overrightarrow{\exp}}
\DeclareMathOperator{\cL}{\mathrsfs{L}}
\DeclareMathOperator{\diag}{diag} \DeclareMathOperator{\e}{e}
\DeclareMathOperator{\id}{Id} 
\DeclareMathOperator{\z}{\mathbf{0}_2}
\newtheorem{fct}{Assumption}
\newtheorem{dfn}{Definition}
\newtheorem{prp}{Proposition}
\newtheorem{thm}{Theorem}[section]

\theoremstyle{remark}

\newcommand{\lo}{L_2(\Omega)}
\newcommand{\W}{\mathring{W}_2^1}
\renewcommand{\endproof}{$\,\,\blacksquare$}
\begin{document}
\selectlanguage{english}
\articolo[Stability of elastic systems]{On stability and
stabilization of elastic systems by time-variant feedback}{M.I.
Caiado\footnotemark[1]\and A.V.
Sarychev\footnotemark[2]\footnotemark[0]} \footnotetext[1]{Centro de
Matem\'{a}tica, Universidade do Minho,
Potugal}\footnotetext[2]{Dipartimento di Matematica per le
Decisioni, Universit\'{a} di Firenze, Italy} \footnotetext[0] {Date:
\today }

\begin{abstract} We study a class of elastic systems
described by a (hyperbolic) partial differential equation. Our
working example is the equation of a vibrating string subject to
linear disturbance. The main goal is to establish conditions for
stabilization and asymptotic stabilization by applying a fast
oscillating control to the  string. In the first situation studied
we assume that system is subject to a damping force; next we
consider the system without damping. We extend the tools of
high-order averaging and of chronological calculus for studying
stability of this distributed parameter system.
\end{abstract}
\section{Introduction}
Stability and stabilization by imposing fast oscillations onto a
system have been extensively studied in literature. Basic example is
stabilization of equilibrium of inverted pendulum by means of fast
harmonic oscillation of its suspension point, \cite{Arn89}. Methods
of high-order averaging developed in \cite{Sary01} are applicable to
study stability and asymptotic stability for a wider class of
(non-harmonically) time-variant systems. As an illustration, a
stabilization condition for pendulum under non-harmonic (fast)
oscillation of its suspension point has been established. An
extension of this work onto the case of planar and spherical double
inverted pendulum was done in \cite{MicS04}.

In the present work we apply the latter approach to a class of
distributed parameter systems: we study  a possibility to stabilize
an elastic system (a string) destabilized by a linear disturbance. A
condition for asymptotic stability under application of a state
feedback control to the damped system is presented. We also obtain a
stability condition for the undamped vibrating string under an
output feedback control.

Our main results are Theorem \ref{Thm:Main} in Section
\ref{sec:String_damp} and Theorem \ref{Thm:Main_undamp} in Section
\ref{sec:String_undamp}. Our problem setting in Section
\ref{sec:String_damp} is close to the one of \cite{kole03} but the
class of stabilizers and the results differ; we provide some comments
on it in Section \ref{Sec:FinalRmk}.

\section{Notation}\label{Sec:Notacao}
Let $\Omega\subset\R$ be a closed interval. As usual $\lo$ is the
space of all measurable (in the sense of Lesbegue) functions on
$\Omega$ having finite norm
\begin{equation*}
    \|u\|= \left(\int_\Omega u^2\, dx\right)^{1/2}.
\end{equation*}
The inner product in $\lo$ is denoted by $<\cdot,\cdot>$. The Hilbert
space consisting of the elements $\lo$ which possess generalized
derivatives up to order $m$ also in $\lo$, is denoted by
$W_2^m(\Omega)$. The norm in this space is denoted by $\|\cdot\|_m$.
Bellow, $\W(\Omega)$ stands for a subspace of the space
$W_2^1(\Omega)$ which coincides with the closure of the set of
infinitely differentiable functions with compact supports in
$\Omega$.

Let $\{\phi_n(x)\}$ be a orthonormal basis of $\lo$. We denote by
$\mathcal{G}_N$ be the subspace of $\lo$ spanned by
$\{\phi_n(x)\}_{n=1}^N$ and by $\mathcal{G}_N^\bot$ the orthogonal
complement of $\mathcal{G}_N$. For an element $v\in \lo$ we denote by
$v_N$ its projection onto $\mathcal{G}_N$ and by $v^N$ the projection
onto $\mathcal{G}_N^\bot$.

\section{Problem setting}\label{Sec:Setting}
We consider forced elastic system described by a (hyperbolic)
partial differential equation. In particular, we consider the
equation for an elastic string
\begin{equation}\label{eq:elastic}
    u_{\tau\tau}=a^2u_{xx}, \qquad x\,\in \,\Omega
\end{equation}
where $u_\xi$ denotes the partial derivative of $u$ with respect to a
variable $\xi$ and the domain $\Omega=[0,2\,\pi]$. We fix boundary
conditions
\begin{equation}\label{eq:boundCond}
    u(x,\tau)=0, \quad x\in \partial \Omega.
\end{equation}

We introduce into \eqref{eq:elastic} a disturbance $\gamma^2\,u$ and
the equation becomes
\begin{equation}\label{eq:string_unst}
     u_{\tau\tau}=a^2u_{xx}+\gamma^2\,u
\end{equation}
whose zero solution is unstable if $\gamma^2>a^2/4$.

If in addition, we assume that there is a damping, then
\eqref{eq:string_unst} turns into
\begin{equation}\label{eq:string}
    u_{\tau\tau}=a^2u_{xx}-\alpha\,u_\tau+\gamma^2\,u, \qquad x\,\in \,\Omega
\end{equation}
with damping coefficient $\alpha>0$.

We would like to achieve stabilization (asymptotic stabilization) of
the null solution of \eqref{eq:string_unst} (\eqref{eq:string}) by
applying time-variant high-gain linear feedback control term
\begin{equation}\label{eq:feed}
 h(u,k\,\tau)=\delta\,k^2\,g(k\,\tau)\,u.
\end{equation}
The resulting ``controlled" equations is
\begin{equation}\label{eq:string_unst_controlled}
    u_{\tau\tau}= a^2 u_{xx}+\gamma^2\,u+h(u,k\,\tau).
\end{equation}
and
\begin{equation}\label{eq:vib_equ}
    u_{\tau\tau}= a^2 u_{xx}-\alpha\,u_\tau+\gamma^2\, u+h(u,k\,\tau) .
\end{equation}
The function $g(\tau)$ is assumed to be bounded, continuous and
$1$-periodic, $\delta>0$ is a small parameter and $k\geq 1$ can be
chosen arbitrarily. Besides, we will assume for $g$:
\begin{fct}\label{assump:1}
$\int_0^1 g(\xi)\,d\xi=0$;
\end{fct}
\begin{fct}\label{assump:2}
For $G(t)=\int_0^t g(\xi)\,d\xi$ there holds
\begin{equation*}
    \int_0^1 G(\xi)\,d\xi=0.
\end{equation*}
\end{fct}

Observe that, under these assumptions,  the time-average of the
vibrational term vanishes, in a way that the averaged equations
\eqref{eq:string_unst_controlled} and \eqref{eq:vib_equ}, coincide
with the equations \eqref{eq:string_unst} and \eqref{eq:string},
respectively. Therefore these averages are unstable.

For given initial data
\begin{equation}\label{eq:inicond}
    u(x,0)=\varphi(x),\qquad u_\tau(x,0)=\psi(x),
\end{equation}
where $\varphi \in W_2^2(\Omega)\cap \W(\Omega)$ and $\psi
\in\W(\Omega)$ the corresponding (classical) solution for
\eqref{eq:string_unst_controlled},\eqref{eq:boundCond} and
\eqref{eq:vib_equ},\eqref{eq:boundCond} exists and is unique in
$W_2^2(Q_T), \, Q_T=\Omega\times(0,T)$, \cite{Lady85}. By a
(classical) solution we mean an element of $W_2^2(\Omega)$ which
satisfies the equation for all $t$ and for almost all $x\in \Omega$.
Next we precise the meaning of stability and asymptotic stability of
the null solution.
\begin{dfn}[Stability]
    The null solution $\tilde{u}(x,t)\equiv0$ of
    \eqref{eq:string_unst_controlled},\eqref{eq:boundCond}
     is stable if for every $\eta>0$ there exists $0<\xi<\eta$ such that
     for any solution $u(x,t)$ with $ \|u(x,0)\|_1^2+\|u_t(x,0)\|^2<\xi$ implies
     \begin{equation*}
    \|u(x,t)\|_1^2+\|u_t(x,t)\|^2<\eta
\end{equation*}
    for all $t\geq 0$ .\endproof
\end{dfn}
\begin{dfn}[Asymptotic stability]
    The null solution $\tilde{u}(x,t)\equiv0$ of (\ref{eq:vib_equ}),(\ref{eq:boundCond})
     is exponentially asymptotically
    stable if it is stable and in addition  there exist
    $\sigma,\,C> 0$ such that
    \begin{equation*}
    \|u(\cdot,t)\|_1^2+\|u_t(\cdot,t)\|^2 \leq C \e^{-\sigma\,t}(\|u(x,0)\|_1^2+\|u_t(x,0)\|^2),
    \end{equation*}
     for all $t\geq 0$.\endproof
\end{dfn}

For convenience we perform the change of the time variable $t=
k\,\tau$. The transformed equation \eqref{eq:string_unst_controlled}
is transformed into
\begin{equation}\label{IBVP_0}
    \begin{split}
&u_{tt}= k^{-2}\,a^2\,u_{xx}+k^{-2}\,\gamma^2\,u+k^{-2}\,h(u,t) \\
&    u(0,t)=u(2\,\pi,t)=0,
\end{split}
\end{equation}
while for the damped equation we obtain
\begin{equation}\label{IBVP_h}
    \begin{split}
&u_{tt}= k^{-2}\,a^2\,u_{xx}-k^{-1}\,\alpha\,u_t+k^{-2}\,\gamma^2\,u+k^{-2}\,h(u,t), \\
&    u(0,t)=u(2\,\pi,t)=0.
\end{split}
\end{equation}
The initial data \eqref{eq:inicond} is converted into
\begin{equation}\label{eq:inicond_2}
  u(x,0)=\varphi(x),\qquad u_t(x,0)=k^{-1}\,\psi(x).
\end{equation}

\section{Chronological calculus: variational formula, logarithm and stability}\label{Sec:ChroCal}
Chronological calculus has been developed by Agrachev and
Gamkrelidze in \cite{AgrGam78} for studying nonlinear time-variant
systems. The systems, to be considered in this work, are linear and
we present some reformulations for the (time-variant) linear case.

Consider the time-variant linear system of ordinary differential
equations
\begin {equation}\label{eq:ODE}
  \dot{z}(t)=A(t)\,z(t),\qquad z(0)=z_0
\end{equation}
where $z(t) \in \mathbb{R}^n$ and $t\mapsto A(t)\in \R^{n\times n}$
is a matrix-valued function depending continuously on $t$.
Considering its solutions $z(t;z_0)$ we can introduce a flow of
linear maps $P^t:z_0\mapsto z(t,z_0)$, or a fundamental matrix
solution. Obviously, $P^t$ is the unique solution to the matrix
ordinary differential equation with initial condition
\begin{equation*}
 \dot{P}^t=A(t)\, P^t,\qquad P^0=\id.
\end{equation*}

Following \cite{AgrGam78}, we call the flow $P^t$ left chronological
exponential of $A(t)$ and denote it
\begin{eqnarray*}
 P^t =\chro \int_0^t A(\tau)\,d\tau.
\end{eqnarray*}
We can also define right chronological exponential
$\overrightarrow{\exp}\int_0^t A(\tau)\,d\tau$ as solution of
$\dot{P}^t=P^tA(t),\, P^0=\id$.

Using this notation, the solution of (\ref{eq:ODE}) is represented as
\begin{equation*}
z(t)= \chro \int_0^t A(\tau)\,d\tau\, z_0=P^t\,z_0.
  \end{equation*}
The flow $\chro \int_0^t A(\tau)\,d\tau$ admits a Volterra series expansion
of the form
\begin{multline*}
     \chro \int_0^t A(\tau)\,d\tau = \id+\int_0^t A(\tau_1)\,d\tau_1+\\
    +\sum_{n=1}^\infty
  \int_0^t\int_0^{\tau_1}\!\!\!\!\ldots\int_0^{\tau_n}
  A(\tau_1) \ldots  A(\tau_n)\,d{\tau_n}\ldots
  d{\tau_1}.
\end{multline*}

\subsection{Variational formula}
Assume now that, in (\ref{eq:ODE}), we take $A(t)= B(t)+ C(t)$,
considering $B(t)$ as a reference matrix of coefficients and $C(t)$
as its perturbation. Then the following chronological calculus
variational formula holds, \cite{AgrGam78},
\begin{multline}\label{eq:VForm}
  \chro \int_0^t [B(\tau)+C(\tau)]\,d\tau =\chro \int_0^t B(\tau)\,d\tau\circ\\\circ\chro \int_0^t \left(\overrightarrow{\exp} \int_0^{\tau} \ad B(\theta)\,d\theta \right)
  C(\tau)\,d\tau.
\end{multline}
Here the operator ``$\ad$'' corresponds to the matrix commutator $
\ad A(t)\,B(t)=- [A(t), B(t)]=B(t)\, A(t)-A(t)\, B(t) $ and
$Q^t=\overrightarrow{\exp} \int_0^t \ad B(\theta)\,d\theta$ is the
solution of the operator differential equation
\begin{equation*}
    \frac{d}{dt}\,Q^t=Q^t \circ \ad B(t), \qquad Q^0=Id.
\end{equation*}

\subsection{Formal expansion for $\ln \chro \int_{0}^{t}
A_{\tau}d\tau$}\label{Sec:Logarithm}
Let $A(t)\stackrel{def}{=} A_t$ be a matrix-valued function with
time-variant entries.  We are interested in a formal expansion for
the logarithm
\begin{equation*}
\Lambda_{0,\,t}(A_\tau)=\ln \chro \int_{0}^{t} A_{\tau}\, d_\tau.
\end{equation*}

It has been shown in \cite{AgrGam78} that $\Lambda_{0,\,t}(A_\tau)$
admits a representation
\begin{equation}\label{eq:Lambda_asympt}
  \Lambda_{0,\,t}(A_\tau)=\sum_{m=1}^{\infty} \Lambda^{(m)}
\end{equation}
where
\begin{equation}\label{eq:Lambda m}
          \Lambda^{(m)}=\int_0^t \!\!\!d\tau_1 \int_0^{\tau_1}\!\!\!\!\!\! d\tau_2 \dots
               \int_0^{\tau_{m-1}}d\tau_m \,p_m(A_{\tau_1}, \dots
               A_{\tau_m}),
\end{equation}
and, for each $m \geq 2$, $p_m(A_{\tau_1}, \dots, A_{\tau_m})$ is a
homogeneous polynomial of first degree in each $A_{\tau_i}$.
Moreover, it is a commutator polynomial in $A_{\tau_1}, \dots,
A_{\tau_m}$, i.e., $p_m$ can be expressed as a linear combination of
$A_{\tau_1}, \dots, A_{\tau_m}$ and of their iterated commutators.

Expressions for the first $\Lambda^{(m)}$ can be founded e.g. in
\cite{Sary01}; in particular $\Lambda^{(1)}$ coincides with the
averaging of $A_\tau$:
\begin{equation}\label{eq:Lambda_1}
    \Lambda^{(1)}=\int_0^t A_\tau\,d\tau.
\end{equation}
 The series
\eqref{eq:Lambda_asympt} is known to be absolutely convergent if
$\int_0^t \|A_\tau\| \,d\tau\leq0.44$, \cite{AgrGam78}.

\subsection{The logarithm and stability of time-variant systems}
Consider linear fast-oscillating system
\begin{align}\label{eq:FastOsc}
 \dot{z}(t)=A(k\,t)\,z(t)
\end{align}
where $x\in\R^n$, $A(\tau)$, $\tau\geq 0,$ is matrix-valued function,
which is  continuous and 1-periodic with respect to $\tau$; $k>0$ is
a large parameter.

A condition for stability of time-periodic system \eqref{eq:FastOsc}
uses the monodromy matrix $P^1$. If all the eigenvalues of $P^1$ are
located in the interior of the unit circle, then the system is
asymptotically stable; system is unstable if at least one eigenvalue
lies outside unit the circle. In general it is difficult to compute
spectrum of $P^1$.

Due to the asymptotic expansion \eqref{eq:Lambda_asympt} it is more
convenient to deal with the logarithm $\ln P^1$. One can arrive to a
conclusion on the system's stability by analyzing the truncations of
the logarithm series if  one gets, somehow, an estimate for the rest
term of such truncation.

In terms of logarithm, stability conditions can be formulated as
follows: system (\ref{eq:FastOsc}) is asymptotically stable if all
the eigenvalues of $\ln P^1$ are located in the open left complex
half-plane and is unstable if at least one eigenvalue lies in the
open right complex half-plane. System (\ref{eq:FastOsc}) is stable if
all the eigenvalues of $\ln P^1$ have non positive real part and
purely imaginary eigenvalues are distinct.

\section{Vibrating string with damping}\label{sec:String_damp}
We consider the problem of asymptotically stabilizing the null
solution of \eqref{IBVP_h} by means of a state feedback control
\begin{equation*}
    h(u,t)=\delta\,k^2\,g(t)\,u.
\end{equation*}
In other words we study asymptotic stability of the null solution of
the equation
\begin{equation}\label{IBVP}
    \begin{split}
&u_{tt}= k^{-2}\,a^2\,u_{xx}-k^{-1}\,\alpha\,u_t+[k^{-2}\,\gamma^2+\delta\,g(t)]\,u, \\
&    u(0,t)=u(2\,\pi,t)=0.
\end{split}
\end{equation}
\subsection{Infinite dimensional system of ODE}
Following the Fourier method, we look for a solution of \eqref{IBVP}
of the form
\begin{equation*}
    u(x,t)=T(t)\,X(x).
\end{equation*}
We will use \emph{dot} and \emph{prime}, respectively, to denote the
derivative with respect to $t$ and to $x$. Differentiating $u$ and
replacing it and its derivatives in the partial differential equation
we obtain, for $X\neq 0$, $T \neq 0$,
\begin{equation}\label{eq:sep_variab}
    \dfrac{\ddot T+k^{-1}\,\alpha\, \dot T+
    b(t)\,T}{k^{-2}\,a^2\,T}=\dfrac{X''}{X}=\lambda,
\end{equation}
where
\begin{equation}\label{eq:b}
    b(t)=-[k^{-2}\,\gamma^2+\delta\, g(t)].
\end{equation}
From \eqref{eq:sep_variab} it is clear that $\lambda$ must be
constant, that is $X''=\lambda\, X$,  $\lambda\in\R$ and
\begin{equation}\label{eq:T}
    \ddot T+k^{-1}\,\alpha\, \dot T+
[b(t)-\lambda\,k^{-2}\,a^2]\,T=0.
\end{equation}
The problem
\begin{equation}\label{eq:spectProb_X}
    X''=\lambda\, X, \quad X(0)=X(2\,\pi)=0,
\end{equation}
 is a ``spectral problem for the elliptic operator"
\begin{equation}\label{eq:L0}
    \cL_0= \dfrac{\partial^2}{\partial x^2}
\end{equation}
whit boundary conditions derived from the original ones in
\eqref{IBVP}. We conclude that in \eqref{eq:spectProb_X}
\begin{align*}
    \lambda_n=-\mu_n^2=-\frac{n^2}{4}, &&
    \tilde{\phi}_n(x)=\sin(\mu_n\,x), \,\,\,n=1,2,\dots
\end{align*}
are the eigenvalues and eigenfunctions. Instead of
$\tilde{\phi}_n(x)$ we take
\begin{equation}\label{eq:base_L2}
    \phi_n(x)=  \pi^{-1/2}
\tilde{\phi}_n(x);
\end{equation}
now $\{\phi_n\}$ forms an orthonormal basis in $\lo$.

We denote by $T_n$ the solution of \eqref{eq:T} with $\lambda$
replaced by $\lambda_n=-\mu_n^2$. The series
\begin{equation}\label{eq:FourSeries}
    u(x,t)=\sum_{n=1}^\infty T_n(t)\, \phi_n(x)
\end{equation}
satisfies, formally, \eqref{IBVP}. For given initial data
\eqref{eq:inicond_2}, we obtain
\begin{align*}
\varphi(x)=\sum_{n=1}^\infty T_n(0)\, \phi_n(x), && &&
\psi(x)=\sum_{n=1}^\infty k\,\dot{T}_n(0)\, \phi_n(x).
\end{align*}
Since $\{\phi_n\}$ is an orthonormal basis in $\lo$, there holds
\begin{align}\label{eq:T_n0}
    T_n(0)=(\varphi,\phi_n) &&
\text{and} && \dot{T}_n(0)=k^{-1}\,(\psi,\phi_n).
\end{align}
Therefore, if $T_n(0)$ and $\dot{T}_n(0)$ are given by
\eqref{eq:T_n0}, then series \eqref{eq:FourSeries} formally satisfies
\eqref{IBVP} with given initial data.

Substituting $u$, given by \eqref{eq:FourSeries} into the partial
differential equation \eqref{IBVP} and equalizing the terms which
contain $\sin(\mu_n\,x)$ we obtain an infinite system of second order
linear ordinary differential equations
\begin{align}\label{eq:T_n}
    \ddot T_n(t)+k^{-1}\,\alpha\, \dot T_n(t)+
[b(t)+\mu_n^2\,k^{-2}\,a^2]\,T_n(t)=0,
\end{align}
for $n=1,2,\dots$. Formally speaking, \eqref{IBVP} is equivalent to
infinite-dimensional system of equations  \eqref{eq:T_n}.

In the next section we illustrate the essence of our method by
studying asymptotic stability of an equation \eqref{eq:T_n}.

\subsection{Stabilization of \eqref{eq:T_n}}\label{Sec:StabOneEq}
If we take time-averaging of the coefficients of \eqref{eq:T_n},
then, due to Assumption \ref{assump:1},  the average of $b(t)$ equals
$-k^{-2}\,\gamma^2$ and we arrive to the equation
\begin{equation}\label{eq:T_n_av}
     \ddot T_n(t)+k^{-1}\,\alpha\, \dot T_n(t)+
k^{-2}(\mu_n^2\,a^2-\gamma^2)\,T_n(t)=0.
\end{equation}

If destabilizing linear disturbance is sufficiently large, e.g.,
\begin{equation*}
    \gamma^2>a^2\, \mu_n^2
\end{equation*}
then the averaged equation \eqref{eq:T_n_av} is unstable. We prove
though that, even in this case, one can achieve stability by choosing
sufficiently large ``frequency" $k$ of (non-harmonic) oscillation
\eqref{eq:feed}.

To derive an asymptotic stability condition for  a single equation in
\eqref{eq:T_n} we first transform it into a two-dimensional
first-order linear system by introducing a new variable $S_n=k\,\dot
T_n$ (c.f. \cite{MicS04}). We get for $z_n^T=(T_n,S_n)$
\begin{equation}\label{eq:system1}
    \dot z_n(t)=A_n(t)\,z_n(t)
\end{equation}
where
\begin{equation*}
    A_n(t)= \left(%
\begin{array}{cc}
  0 & k^{-1} \\
  k^{-1}\,(\gamma^2-\mu_n^2\, a^2)+k\,\delta\, g(t)& -k^{-1}\,\alpha \\
\end{array}%
\right).
\end{equation*}
We represent $A_n(t)$ as $B_n(t)+C_n$, where
\begin{align*}
    &B_n(t)=\delta\,k\,g(t)\,B, \qquad B=\left(%
\begin{array}{cc}
  0 & 0 \\
  1 & 0 \\
\end{array}%
\right),
\end{align*}
{and}
\begin{align*}
&C_n=k^{-1}\left(%
\begin{array}{cc}
  0 & 1 \\
  \gamma^2-a^2\,\mu_n^2 & -\alpha \\
\end{array}%
\right).
\end{align*}

According to the chronological calculus variational formula
\eqref{eq:VForm}, the monodromy matrix of the system
\eqref{eq:system1} can be represented as
\begin{multline*}
    P^1_n=\chro \int_0^1 A_n(\xi)\,d\xi=\chro \int_0^1 [B_n(\xi)+C_n]\,d\xi\\
    = \chro \int_0^1 B_n(\tau)\,d\tau
    \circ\chro \int_0^1 \left(\overrightarrow{\exp} \int_0^{\tau} \ad
    B_n(\theta)\,d\theta\,
  C_n\right)\,d\tau.
\end{multline*}
From the Assumption \ref{assump:1}, it follows that
\begin{align*}
    \int_0^1
B_n(\tau)\,d\tau=\z && \text{and, hence} && \chro \int_0^1
B_n(\tau)\,d\tau= \id_2,
\end{align*}
where $\z$ and $\id_2$ are, respectively, the two-dimensional null
and identity matrices. Besides
\begin{equation*}
 \overrightarrow{\exp} \left(\int_0^{t} \ad
    B_n(\theta)\,d\theta\right)\,   C_n=\e^{-\delta\,k\, G(t)\ad B}C_n
\end{equation*}
where $\int_0^t g(\theta)\,d\,\theta=G(t)$ (c.f. Assumption
\ref{assump:2}).

We represent the monodromy matrix as
\begin{equation*}
    P^1_n=\chro \int_0^1 D_n(t)\,d\,t,
\end{equation*}
where
\begin{align*}
  D_n(t)&=\e^{-\delta\,k\, G(t)\ad B}C_n\\
    &=[\id_2-\delta\,k\, G(t)\ad B+ \frac{(\delta\,k)^2}{2}\, G^2(t)\ad^2
    B-
    \frac{(\delta\,k)^3}{3!}\,G(t)\ad^3 B+ \dots ]\,C_n.
\end{align*}

As far as $\ad^j\, B\,C_n=\z,\,j\geq 3$, the series for
$\e^{-\delta\,k\, G(t)\ad B}C_n$  ends at the third term and
therefore
\begin{equation}\label{eq:D_n}
\!\!\!D_n(t)\!=\!\left(\!\!\!%
\begin{array}{cc}
  \delta\,G(t) & k^{-1} \\
  k^{-1}(\gamma^2-a^2\mu_n^2)-\tilde G(t) & -k^{-1}\alpha-\delta\, G(t) \\
\end{array}%
\!\!\!\right) \!\!\!\!\!\! \!\!\!\!\!\!
\end{equation}
where
 $\tilde G(t)=\alpha\, \delta\, G(t)+k\,\delta^2\, G^2(t)$.

Let now
\begin{equation*}
    \Lambda=\ln P^1_n= \ln\left(\chro \int_0^1 D_n(t)\,dt\right)
\end{equation*}
be the logarithm of the monodromy matrix.

Consider the expansion \eqref{eq:Lambda_asympt} for the logarithm,
with $\Lambda^{(m)}$ defined by \eqref{eq:Lambda m}. We will show
that, under some conditions, the stability properties of
\eqref{eq:system1} are determined by $\Lambda^{(1)}$.

As we said in Section \ref{Sec:Logarithm}, $\Lambda^{(1)}$ coincides
with the averaging of $D_n(t)$.  According to Assumption
\ref{assump:2}
\begin{align}\label{eq:Lambda_1_damp}
    \Lambda^{(1)}=\left(%
\begin{array}{cc}
   0 & k^{-1} \\
  k^{-1}(\gamma^2-a^2\,\mu_n^2)-\delta^2\,k\,\Gamma & -k^{-1}\alpha \\
\end{array}%
\right)
\end{align}
where
\begin{equation}\label{eq:Gama}
    \Gamma=\int_0^1 G^2(t)\,dt>0.
\end{equation}

The eigenvalues of $\Lambda^{(1)}$ can be either real or complex
conjugated; they are equal to
\begin{equation}\label{eq:eigenVL1}
    \xi^{(1)}_i=\frac{\Theta \pm \sqrt{\Theta^2-4 \Delta}}{2}
\end{equation}
where $\Theta=-k^{-1}\,\alpha$ and
$\Delta=\delta^2\,\Gamma-k^{-2}(\gamma^2-a^2\,\mu_n^2)$ are the
trace and the determinant of $\Lambda^{(1)}$, respectively.
Obviously, $\Theta<0$. $\Lambda^{(1)}$ possess real negative
eigenvalues if
\begin{equation}\label{eq:negatEV real}
    k^{-2}\frac{4\gamma^2-a^2\,n^2}{4}<\delta^2\int_0^1 G^2 (t)\, dt
    \leq k^{-2}\frac{\alpha^2+4\gamma^2-a^2\,n^2}{4}.
\end{equation}
For large $n$ the latter inequality is never satisfied and then for
\begin{equation}\label{eq:negatEV compl}
    \delta^2\int_0^1 G^2 (t)\, dt>
k^{-2}\frac{\alpha^2+4\gamma^2-a^2\,n^2}{4}
\end{equation}
we obtain a pair of conjugate complex numbers with negative real
part.

What for the rest term $\Lambda-\Lambda^{(1)}$, then the following
estimate  is available.
\begin{prp}[\cite{MicS04}]\label{Thm:prop 1}
For $n^2<4\frac{k^2+\gamma^2}{a^2}$ there exists a positive constant
$c$ and a $2\times2$ matrix $R$ such that
\begin{equation}\label{eq:error}
    \Lambda-\Lambda^{(1)}= \delta^2 \,R, \quad \|R\|_M < c
\end{equation}
where neither $c$ nor $R$ depend on $n$.\endproof
\end{prp}

Here $ \|\cdot\|_M$ stands for some matrix norm.  The upper bound on
$n$ is irrelevant; we can assume it to be fulfilled  for $k>k_0$ as
stated in the next proposition. The estimate \eqref{eq:error} can be
derived from convergence results for the chronological exponential
presented in \cite{AgrGam78}; in \cite{MicS04} the detailed
exposition for the particular case of stability of inverse double
pendulum can be found.

From \eqref{eq:negatEV real}, \eqref{eq:negatEV compl} and
\eqref{eq:error}, we arrive to the condition for the asymptotic
stability of  the two-dimensional system \eqref{eq:system1}.

\begin{prp}\label{Thm:NegEV_N}
 For each $\epsilon>0$ there exist $\delta_0>0,\, k_0>1$ such that
 the null solution $z_n(t)\equiv 0$ of \eqref{eq:system1}
 is exponentially asymptotically stable if
$0<\delta<\delta_0, k>k_0$ and
\begin{equation}\label{eq:AspSta z_n}
    \delta^2\int_0^1 G^2 (t)\, dt>
k^{-2}\frac{4\gamma^2-a^2\,n^2}{4}+\epsilon,
\end{equation}
and unstable if $0<\delta<\delta_0, k>k_0$ and
\begin{equation}\label{eq:UnsSta z_n}
    \delta^2\int_0^1 G^2 (t)\, dt<
k^{-2}\frac{4\gamma^2-a^2\,n^2}{4}-\epsilon.\,\,\, \blacksquare
\end{equation}
\end{prp}

Observe that if the condition \eqref{eq:AspSta z_n} is fulfilled for
$n=N$ then it is fulfilled for all $n \geq N$.

\subsection{Asymptotic stability of Galerkin's approximation}\label{Sec:Galerkin}
The Galerkin approximation for \eqref{IBVP} is a system of
(decoupled) equations \eqref{eq:T_n} for $n=1, \ldots,N$.

This system is equivalent to a $2\,N$ block diagonal system
\begin{equation}\label{eq:syst_N}
    \dot z_N(t)= H_N(t)\,z_N(t)
\end{equation}
where $z_N^T(t)=(\tilde{z}_1(t),\dots,\tilde{z}_N(t))$ and
\begin{equation*}
    H_N(t)=\diag \{ D_1(t), \dots, D_N(t)\}.
\end{equation*}
Due to the decoupling, we may conclude that the system
\eqref{eq:syst_N} is asymptotically stable whenever each of the
systems
\begin{equation*}
    \dot z_n(t)=D_n(t)\,z_n(t), \qquad n=1,\dots, N,
\end{equation*}
is asymptotically stable; it is unstable if at least one of them is
unstable. Instability occurs if \eqref{eq:UnsSta z_n} holds for some
$n=n_1\leq N$ (then it holds for all $n\leq n_1$). According to
Proposition \ref{Thm:NegEV_N} system \eqref{eq:syst_N} is
exponentially asymptoticaly stable if
\begin{equation}\label{eq:AsympStab_N}
    \delta^2\int_0^1 G^2 (t)\, dt>
k^{-2}\frac{4\gamma^2-a^2}{4}+\epsilon.
\end{equation}
Observe that the latter condition does not depend on $N$.

Let $\mathcal{G}_N$ be the subspace of $\lo$ spanned by
$\{\phi_n(x)\}_{n=1}^N$ defined in \eqref{eq:base_L2} and let
$\mathcal{G}_N^\bot$ be the orthogonal complement of $\mathcal{G}_N$
in $\lo$.

We denote the projection of the (unique) solution of \eqref{IBVP} on
$\mathcal{G}_N$ by $u_N$ and the projection of $u$ onto
$\mathcal{G}_N^\bot$ by $u^N$.

 Then
\begin{equation*}
    u_N(x,t)=\sum_{n=1}^N (u, \phi_n)\,\phi_n(x)=\sum_{n=1}^N
    T_n(t)\phi_n(x)
\end{equation*}
 is the solution of the $N$th Galerkin's approximation of
\eqref{IBVP} and the following proposition holds.
\begin{prp}\label{Pro:Prop3}
    For each $\epsilon>0$ there exist $\delta_0>0,\, k_0>0$ if
$0<\delta<\delta_0, k>k_0$ and
\begin{equation}\label{eq:AsymSta z_N}
    \delta^2\int_0^1 G^2 (t)\, dt>
k^{-2}\frac{4\gamma^2-a^2}{4}+\epsilon
\end{equation}
then     the null solution  $z_n(t)\equiv 0$ of the Galerkin's
approximation of \eqref{IBVP}
     is exponentially asymptotically stable
\endproof
\end{prp}

Therefore, if \eqref{eq:AsymSta z_N} holds then
\begin{align*}
    &&\|u_N(\cdot,t)\|_1^2+\|(u_N)_t(\cdot,t)\|^2\rightarrow 0
    &&\text{as}&&t\rightarrow +\infty.&&
\end{align*}
 To establish asymptotic stability of the
problem \eqref{eq:vib_equ},\eqref{eq:boundCond} it suffices to
prove, that for sufficiently large $N$, the  projection of $u$ onto
$\mathcal{G}^\bot_N$ tends to zero as $t\rightarrow +\infty$
\begin{align*}
  &&  \|u^N(\cdot,t)\|_1^2+\|u^N_t(\cdot,t)\|^2\rightarrow 0 && \text{as}&& t\rightarrow
    +\infty.&&
\end{align*}

\subsection{Convergence of $\|u^N(\cdot,t)\|_1^2+\|u^N_t(\cdot,t)\|^2$}\label{Sec:GalerkinConv}
As far as $u$ is an element of $W_2^2(Q_T)$, where
$Q_T=\Omega\times(0,T)$, $u^N(\cdot,t)$ is in $W_2^2(\Omega)$ for all
$t\in(0,T)$. We prove now
\begin{prp}\label{Pro:Prop4}
For the projection $u^N$ there exists $c>0$ such that
\begin{equation*}
    \|u^N(\cdot,t)\|_1^2+\|u^N_t(\cdot,t)\|^2\leq \e^{-c\,t}(\|u^N(\cdot,0)\|_1^2+\|u^N_t(\cdot,0)\|^2),
\end{equation*}
for all $t\geq0$.\endproof
\end{prp}

\proof Introduce  a function $V$ defined, for any $u\in W_2^2(Q_T)$,
by
\begin{equation}\label{eq:Liap_damp}
    V=\frac{k^{-2}\alpha^2}{4}\|u\|^2+
    \frac{k^{-2}a^2}{2}\|u_x\|^2
    +\frac{1}{2}\|u_t\|^2+\frac{k^{-1}\alpha}{2}<u,u_t>.
\end{equation}
We recall that $<\cdot,\cdot>$ and $\|\cdot\|$  stand for the inner
product and for the norm in $\lo$, respectively (c.f. Section
\ref{Sec:Notacao}). Let $V_N$ denote the restriction of the function
$V$ onto $\mathcal{G}_N^\bot$.

We first note that there is a positive constant $C_1$ such that
\begin{equation}\label{eq:C1}
 V_N\geq C_1(\|u^N\|_1^2+\|u^N_t\|^2).
\end{equation}
 Indeed, $V$ can be represented as
\begin{equation*}
    V=\frac{k^{-2}a^2}{2}\|u_x\|^2+ \frac{1}{4}\|u_t\|^2+
    \frac{1}{4}\left\|k^{-1}\alpha\,u+u_t\right\|^2.
\end{equation*}
As long as the inequality $\|u^N_x\|^2\geq N^2 \|u^N\|^2$ holds,
then
\begin{equation*}
     V_N\geq\frac{k^{-2}a^2}{4}N^2\|u^N\|^2+\frac{k^{-2}a^2}{4}\|u_x^N\|^2
     + \frac{1}{4}\|u_t^N\|^2+
    \frac{1}{4}\left\|k^{-1}\alpha\,u^N+u_t^N\right\|^2
\end{equation*}
and \eqref{eq:C1} follows with
$C_1=C_1(k)=\min\left(\frac{k^{-2}a^2}{4}, \frac{1}{4}\right)$.

Let  $\dot V_N$ denote the time derivative of $V_N(u^N(\cdot,t))$.
After some simple computation we conclude
\begin{equation*}
    -\dot V_N
    \geq C_2(N,k,t)\|u^N\|^2+
    \frac{k^{-3}\alpha}{4}\,a^2\|u^N_x\|^2+\frac{k^{-1}\alpha}{4}\|u^N_t\|
\end{equation*}
where
\begin{equation*}
     C_2(N,k,t)=
    \frac{k^{-3}\alpha}{4}\,a^2\,N^2-
    \left[\frac{k^{-1}\alpha}{2}b(t)+\frac{1}{k^{-1}\alpha}\,b^2(t)\right].
\end{equation*}
As long as $g(t)$ and $b(t)$ are periodic and bounded, there exists
$\beta>0$ such that $|b(t)|\leq \beta$ for all $t\geq 0$, (c.f.
\eqref{eq:b}). Then, for all $t\geq 0$
\begin{equation*}
    C_2(N,k,t)\geq \hat{C}_2(N,k)=
    \frac{k^{-3}\alpha}{4}\,a^2\,N^2-
    \left[\frac{k^{-1}\alpha}{2}\beta+\frac{\beta^2}{k^{-1}\alpha}\right].
\end{equation*}
We assume $N$ to be sufficiently large to guarantee
$\hat{C}_2(N,k)>0$ (for all $\delta< \delta_0$, $k>k_0$). Then, for
those $N$, $\hat{C}_2$ can be chosen independently on $N$
\begin{equation*}
     -\dot V_N
    \geq \hat{C}_2(N,k)\|u^N\|^2+
    \frac{k^{-3}\alpha}{4}\,a^2\|u^N_x\|^2+\frac{k^{-1}\alpha}{4}\|u^N_t\|>0.
\end{equation*}
We can also show that
\begin{equation*}
     V_N\leq \frac{3}{8}k^{-2}\alpha^2\|u^N\|^2+
    \frac{k^{-2}}{2}\,a^2\|u^N_x\|^2+\|u^N_t\|.
\end{equation*}
Therefore, there exists a positive constant $C(k)$ such that $\dot
V_N\leq -C(k)\,V_N$, hence
\begin{equation}\label{eq:C}
 V_N(t)\leq V_N(0)\e^{-C(k)\,t}.
\end{equation}
 From \eqref{eq:C1} and \eqref{eq:C} we conclude
\begin{equation*}
    \|u^N(\cdot,t)\|_1^2+\|u^N_t(\cdot,t)\|^2\rightarrow0 \qquad \text{as}\qquad t\rightarrow+\infty
    .\,\,\, \blacksquare
\end{equation*}

From  the Propositions \ref{Pro:Prop3} and \ref{Pro:Prop4}  we
conclude the asymptotic stability of the null solution of infinite
dimensional system \eqref{eq:syst_N} and the consequent asymptotic
stability of the null solution of the problem
(\ref{eq:vib_equ}),(\ref{eq:boundCond}).

\subsection{Main result for damped string}
Sumarizing the reasoning of the previous sections  we state

\begin{thm}\label{Thm:Main}
Consider the equations  \eqref{eq:vib_equ},\eqref{eq:boundCond}
 for damped elastic string subject to linear
disturbance and controlled by a time-variant (high-gain
fast-oscillating non-harmonic) linear feedback control
\eqref{eq:feed}.

    For each $\epsilon>0$ there exist $\delta_0>0,\, k_0>0$ such that if
$0<\delta<\delta_0, k>k_0$ and the inequality \eqref{eq:AsymSta z_N}
holds, then   the null solution $u(x,t)\equiv0$ of problem
    \eqref{eq:vib_equ},\eqref{eq:boundCond}
     is exponentially asymptotically stable. \endproof
\end{thm}


\section{Undamped vibrating string}\label{sec:String_undamp}
In this section, we consider the problem of stabilization of the
null solution of undamped equation \eqref{eq:string_unst}. We try
to stabilize it by a forcing term, which results in equation
\begin{equation}\label{IBVP_0_out}
    \begin{split}
&u_{tt}= k^{-2}\,a^2\,u_{xx}+k^{-2}\,\gamma^2\,u+k^{-2}\,h(u,t) \\
&    u(0,t)=u(2\,\pi,t)=0.
\end{split}
\end{equation}

Our choice is a time-variant ``output" feedback, \cite{Cor99},
\begin{equation}\label{eq:feed_out}
    h(u,t)=\delta\,k^2\,g(k\,\tau)\,\Pi_N\,u,
\end{equation}
 c.f. \eqref{eq:feed}, where $\Pi_N$ stands for the orthogonal
projection of the ``state" $u$ onto the subspace of $\lo$
generated by the first $N$ harmonics.

For the control $h$ defined by \eqref{eq:feed_out}, the equation
 \eqref{IBVP_0_out} is formally equivalent to an infinite
dimensional system of second order linear ordinary differential
equations
\begin{align}
  &\ddot T_n(t)+[k^{-2}\,(\mu_n^2\,a^2-\gamma^2)-\delta\,g(t)]\,T_n(t)=0,&&  n\leq N \label{eq:T_n_undam}\\
    &\ddot T_n(t)+k^{-2}\,(\mu_n^2\,a^2-\gamma^2)\,T_n(t)=0,&&
    n> N.  \nonumber
\end{align}
This system is obtained by introducing the (formal) expansion
\eqref{eq:FourSeries} of $u$.

The analysis of the stability of $u_N=\Pi_N u$ goes along with
what was done in Sections \ref{Sec:StabOneEq} and
\ref{Sec:Galerkin} for the damped case.

As before, we introduce a new variable $S_n=k\,\dot T_n$ and reduce
\eqref{eq:T_n_undam} to a two-dimensional first-order linear system
\begin{equation}\label{eq:system_n_umd}
    \dot z_n(t)=\tilde A_n(t)\,z_n(t)
\end{equation}
where
\begin{equation*}
    \tilde A_n(t)= \left(%
\begin{array}{cc}
  0 & k^{-1} \\
  k^{-1}\,(\gamma^2-\mu_n^2\, a^2)+k\,\delta\, g(t)& 0\\
\end{array}%
\right).
\end{equation*}
Using chronological calculus variational formula \eqref{eq:VForm}, we
compute the monodromy matrix for \eqref{eq:system_n_umd}
\begin{equation*}
    Q^1_n=\chro \int_0^1 \tilde D_n(t)\,d\,t,
\end{equation*}
where
\begin{equation*}
\tilde D_n(t) =\left(
\begin{array}{cc}
  \delta\,G(t) & k^{-1} \\
  k^{-1}(\gamma^2-a^2\mu_n^2)-k\,\delta^2\, G^2(t) & -\delta\, G(t) \\
\end{array}\right)
\end{equation*}
and $G(t)=\int_0^t g(\xi)\,d\xi$. Next we consider the logarithm
$\Lambda=\ln Q^1_n$ given by the asymptotic expansion
\eqref{eq:Lambda_asympt} and compute $\Lambda^{(1)}$
\begin{align*}
   && \Lambda^{(1)}=\left(%
\begin{array}{cc}
   0 & k^{-1} \\
  k^{-1}(\gamma^2-a^2\,\mu_n^2)-\delta^2\,k\,\Gamma & 0 \\
\end{array}%
\right),
\end{align*}
where $\Gamma=\int_0^1 G^2(t)\,dt>0$. It is clear that the main
difference between damped and undamped systems is that now
$\Lambda^{(1)}$ is a traceless matrix while
\eqref{eq:Lambda_1_damp} had negative trace. Hence, to achieve
stability we have to ensure that its eigenvalues are two imaginary
(conjugate) numbers. Stability condition for (each of) the
two-dimensional system \eqref{eq:system_n_umd} is analogous to the
one of the Proposition \ref{Thm:NegEV_N}.
\begin{prp}\label{Thm:NegEV_N_undamp}
 For each $\epsilon>0$ there exist $\delta_0>0,\, k_0>1$ such that
 the null solution $z_n(t)\equiv 0$ of \eqref{eq:system_n_umd}
 is exponentially asymptotically stable if
$0<\delta<\delta_0, k>k_0$ and
\begin{equation*}\label{eq:AsymSta z_n_umd}
    \delta^2\int_0^1 G^2 (t)\, dt>
k^{-2}\frac{4\gamma^2-a^2\,n^2}{4}+\epsilon
\end{equation*}
and unstable if $0<\delta<\delta_0, k>k_0$ and
\begin{equation*}\label{eq:UnsSta z_n_umd}
    \delta^2\int_0^1 G^2 (t)\, dt<
k^{-2}\frac{4\gamma^2-a^2\,n^2}{4}-\epsilon.\,\,\, \blacksquare
\end{equation*}
\end{prp}
Recall that Galerkin's approximation for the equation
\eqref{IBVP_0_out} is a system of $N$ decoupled second order
equations \eqref{eq:system_n_umd}, where  $n=1, \cdots N$. As in
Section \ref{Sec:Galerkin}, this system is stable whenever all the
equations are stable, that is, when the Proposition
\ref{Thm:NegEV_N_undamp} holds for all $n \leq N$.

To be able to conclude stability of the null solution of
\eqref{IBVP_0_out}, \eqref{eq:feed_out} we have to analyze stability
of  the projection of \eqref{IBVP_0_out} onto $\mathcal{G}_N^\bot$
\begin{equation*}
    u_{tt}^N= k^{-2}\,a^2\,u_{xx}^N+k^{-2}\,\gamma^2\,u^N.
\end{equation*}

Let us consider a function
\begin{equation*}
    V=-k^{-2}\frac{\gamma^2}{2}\|u\|^2+k^{-2}\frac{a^2}{2}\|u_x\|^2+\frac{1}{2}\|u_t\|^2.
\end{equation*}
This function is not necessarily positive but its restriction
$V_N$ onto $\mathcal{G}_N^\bot$ is positive provided that $N$ is
sufficiently large, say,  $N^2> 2\gamma^2/a^2$. Indeed, since the
inequality $\|u_x^N\|^2\geq N^2\|u^N\|^2$ holds, then in this case
\begin{equation*}
   V_N\geq \frac{k^{-2}}{2}\left[a^2\,N^2-\gamma^2\right]\|u^N\|^2+
    k^{-2}\frac{a^2}{4}\|u_x\|^2+\frac{1}{2}\|u_t^N\|^2\geq    0.
\end{equation*}
Let $\dot V_N$ denote the time derivative of $V_N(u^N(\cdot,t))$.
Then
\begin{align*}
    \dot
    V_N&=-k^{-2}\,\gamma^2<u^N,u_t^N>+k^{-2}a^2<u_x^N,u_{xt}^N>+<u_t^N,u_{tt}^N>\\
    &=-k^{-2}\,\gamma^2<u^N,u_t^N>-k^{-2}a^2<u_{xx}^N,u_{t}^N> +k^{-2}\,a^2\,
    <u_t^N,u_{xx}^N>+\\
    &\qquad+k^{-2}\,\gamma^2\,<u_t^N,u^N>\\
    &=0.
\end{align*}
The equality $<u_x^N,u_{xt}^N>=-<u_{xx}^N,u_{t}^N>,$ which holds
under the boundary conditions \eqref{eq:boundCond} is concluded by
integration by parts with respect to $x$.

Therefore $V_N(u^N(\cdot,t)) \equiv \mbox{const}$ wherefrom we
 obtain an upper estimate for $u^N$ and conclude stability of the null
solution of of  the projection of \eqref{IBVP_0_out} onto
$\mathcal{G}_N^\bot$. Completing this by the result of the
Proposition~\ref{Thm:NegEV_N_undamp} we conclude stability of
\eqref{IBVP_0_out},\eqref{eq:boundCond}.

\begin{thm}\label{Thm:Main_undamp}
Consider the equations
\eqref{eq:string_unst_controlled},\eqref{eq:boundCond}
 for an undamped elastic string which is subject to linear
disturbance and is controlled by a time-variant (high-gain
fast-oscillating non-harmonic) linear output feedback control
\eqref{eq:feed_out}. Then for each $\epsilon>0$ there exist
$\delta_0>0,\, k_0>0$ such that if $0<\delta<\delta_0, k>k_0$ and the
inequality
\begin{equation*}
    \delta^2\int_0^1 G^2 (t)\, dt>
k^{-2}\frac{4\gamma^2-a^2}{4}+\epsilon
\end{equation*}
holds, then   the null solution  of problem
\eqref{eq:string_unst_controlled},\eqref{eq:boundCond} is stable.
\endproof
\end{thm}

\section{Final remarks}\label{Sec:FinalRmk}
From \eqref{eq:eigenVL1}, it is clear that the degree of asymptotic
stability increases with $\alpha$. An interesting but not unexpected
phenomena arises when $\Theta^2-4 \Delta$ equals zero: multiple
negative eigenvalue bifurcates to a conjugate complex pair.

If \eqref{eq:UnsSta z_n} holds for some $n=n_1 \leq N$, then finite
dimensional approximation \eqref{eq:syst_N} is unstable implying that
the null solution of the problem
(\ref{eq:vib_equ}),(\ref{eq:boundCond}) is unstable.

The output feedback control constructed in the
Section~\ref{sec:String_undamp} is suitable for the asymptotic
stabilization of the damped string (system
(\ref{eq:vib_equ}),(\ref{eq:boundCond})).

The approach to the problem presented in this work is  more general
than the one of \cite{kole03}. The main difference is that we do not
select a particular type of vibrational control, in contrast to the
harmonic vibration and rather restrictive assumption for
frequencies, which appeared in \cite{kole03}. Besides here we were
able to deal with stabilization of an undamped string by means of
output feedback.

\section{Acknowledgments}
A.Sarychev has been partially supported by MIUR, Italy, COFIN
program, grant 2004015409-003. M.I. Caiado was partially supported
by Marie Curie Control Training Site (HPMT-CT-2001-00278-82) and by
a short-term fellowship by Funda\c{c}\~{a}o Calouste Gulbenkian (Portugal).
She expresses her gratitude to the Dipartimento di Matematica per le
Decisioni, Universit\`a di Firenze (Italy), for the hospitality
offered during her visits to the department.

\bibliographystyle{ieeetran}
\bibliography{CaiSar}
\bigskip
\begin{flushleft}

{\bf AMS Subject Classification: 34E05, 34C29, 34D20.}\\[2ex]

A.V. SARYCHEV\\
Dipartimento di Matematica per le Dicisione\\
Universit\`{a} di Firenze\\
Via C. Lombroso, 6/7\\
50134  Firenze, Italia\\
e-mail: \texttt{andrey.sarychev@dmd.unifi.it}\\[2ex]

M.I. CAIADO\\
 {Centro de Matem\'{a}tica}\\
 {Universidade do Minho}\\
{Campus de Gualtar} \\
4710-057 Braga,  Portugal\\
e-mail: \texttt{icaiado@math.uminho.pt}\\[2ex]

\end{flushleft}
\end{document}